\documentclass{article}

\usepackage{arxiv}

\usepackage[utf8]{inputenc} % allow utf-8 input
\usepackage[T1]{fontenc}    % use 8-bit T1 fonts
\usepackage{hyperref}       % hyperlinks
\usepackage{url}            % simple URL typesetting
\usepackage{booktabs}       % professional-quality tables
\usepackage{amsfonts}       % blackboard math symbols
\usepackage{nicefrac}       % compact symbols for 1/2, etc.
\usepackage{microtype}      % microtypography
\usepackage{graphicx}
\usepackage[numbers]{natbib}
\usepackage{doi}

\usepackage{mathtools}
\usepackage{amsthm}

\usepackage{multirow}
\usepackage{csquotes}
	
\usepackage{lscape}

\title{Energy System Optimisation using (Mixed Integer) Linear Programming}

\author{ {Sebastian~Miehling}\thanks{Dual first authorship} \\
	Chair of Energy Systems \\
	TUM School of Engineering and Design\\
	Technical University of Munich\\
	\And
	\href{https://orcid.org/0000-0003-0135-6167}{\includegraphics[scale=0.06]{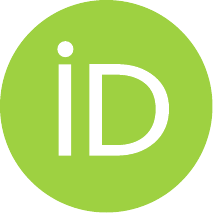}}\hspace{1mm}{Andreas~Hanel}$^*$ \\
	Chair of Energy Systems\\
	TUM School of Engineering and Design\\
	Technical University of Munich\\
	\texttt{andreas.hanel@tum.de}
	\And
	\href{https://orcid.org/0000-0002-2028-6769}{\includegraphics[scale=0.06]{orcid.pdf}}\hspace{1mm}{Jerry~Lambert} \\
	Chair of Energy Systems\\
	TUM School of Engineering and Design\\
	Technical University of Munich\\
	\And
	\href{https://orcid.org/0000-0002-9100-7456}{\includegraphics[scale=0.06]{orcid.pdf}}\hspace{1mm}{Sebastian~Fendt} \\
	Chair of Energy Systems\\
	TUM School of Engineering and Design\\
	Technical University of Munich\\
	\texttt{sebastian.fendt@tum.de}
	\And
	\href{https://orcid.org/0000-0003-3743-2731}{\includegraphics[scale=0.06]{orcid.pdf}}\hspace{1mm}{Hartmut~Spliethoff} \\
	Chair of Energy Systems\\
	TUM School of Engineering and Design\\
	Technical University of Munich\\
}

\renewcommand{\shorttitle}{Optimal Technology Utilisation in Multi-sectoral Applications}

\hypersetup{
pdftitle={OpTUMus},
pdfauthor={M.~Miehling, A.~Hanel, J.~Lambert, F.~Fendt, H.~Spliethoff}
}

\begin{document}
\maketitle

\begin{abstract}
Although energy system optimisation based on linear optimisation is often used for influential energy outlooks and studies for political decision-makers, the underlying background still needs to be described in the scientific literature in a concise and general form. This study presents the main equations and advanced ideas and explains further possibilities mixed integer linear programming offers in energy system optimisation. Furthermore, the equations are shown using an example system to present a more practical point of view. Therefore, this study is aimed at researchers trying to understand the background of studies using energy system optimisation and researchers building their implementation into a new framework. This study describes how to build a standard model, how to implement advanced equations using linear programming, and how to implement advanced equations using mixed integer linear programming, as well as shows a small exemplary system.
\begin{itemize}
    \item Presentation of the OpTUMus energy system optimisation framework
    \item Set of equations for a fully functional energy system model
    \item Example of a simple energy system model
\end{itemize}
\end{abstract}

% keywords can be removed

\section{Introduction}
Linear programming for modelling energy systems is widespread and state-of-the-art. The goal is to investigate the impact of various boundary conditions and technological developments to achieve a sustainable energy supply. A well-known model is the world energy model (WEM), which is used by the international energy agency (IEA), for example, as the basis for the world energy outlook \cite{Hughes.2021}. The WEM is a large-scale model of the global energy system that can be used to analyse the individual sectors in different regions of the world. Another tool for evaluating energy systems is TIMES, developed by an international consortium within the framework of IEA-ETSAP \cite{Loulou.2016}. TIMES uses the approach of minimising the total costs under consideration of various constraints, e.g., technologies, political goals, or physical laws. The tool is used, among others, for scenario-based analyses of possible energy system developments. Fraunhofer ISE developed REMod to investigate possible developments of national energy systems \cite{Erlach.2018}. The model is equally based on the minimisation of total costs, whereby the optimisation problem may have non-linear properties.

However, the models are usually presented in the context of explicit problems. Thereby the explanation of the basic methodology often comes too short. One reason is that these models are extremely large and complex, and more space is needed in the actual publications. Nevertheless, or even for this reason, readers from other disciplines often have problems understanding the methodology in its entirety.

Therefore, this study aims to present the basic concepts of energy system optimisation, the underlying mathematical formulation, and an exemplary system. This study also explains how advanced features can be implemented into this kind of model by using further constraints or mixed integer linear programming. The implementation of the model is based on the framework “Optimal Technology Utilisation in Multi-Sectoral Applications” (OpTUMus). This framework was developed by the Chair of Energy Systems (Technical University of Munich) and used in \cite{Miehling.2022} and \cite{Miehling.2021}. 

% ----------------------------------------------------------------------------------------------------------------------------
% ----------------------------------------------------------------------------------------------------------------------------
% ----------------------------------------------------------------------------------------------------------------------------
% ----------------------------------------------------------------------------------------------------------------------------
% ----------------------------------------------------------------------------------------------------------------------------
\section{Method Details}
The proposed energy system optimisation method aims to find the ideal scheduling of components and, if needed, their installed capacity. Describing the system as a graph using nodes and edges is purposeful. A conservation equation must be fulfilled for each time step and at each node. These conservation equations can refer to energy flows, e.g., electricity or heat, and material flows, such as hydrogen, natural gas, or chemicals.

\begin{figure}[hbtp]
    \centering
    \includegraphics[width=0.5\textwidth]{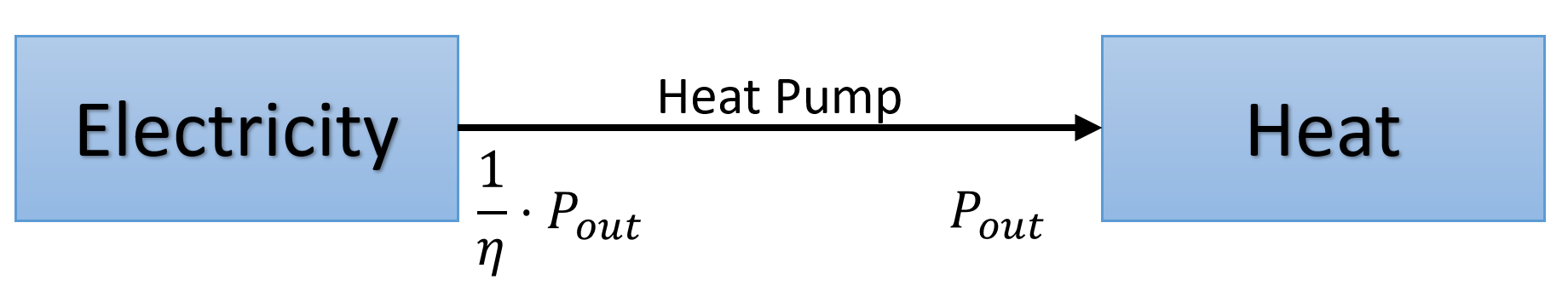}
    \caption{A Heat Pump Converting Electricity to Heat}
    \label{fig:001}
\end{figure}

The edges of the graph represent the components of the energy system. These edges transport an energy or material flow from one node to another or across the system boundary. An example of an edge (component) is a power plant that produces electricity and thus supplies electricity to the corresponding node. Another example is a heat pump, which takes energy from the electricity node and supplies heat to a heat node with an exemplary coefficient of performance of 4. \autoref{fig:001} shows how this heat pump would transfer energy from the electricity node to the heat node. The edges can have more than one starting or endpoint.

\begin{figure}[hbtp]
    \centering
    \includegraphics[width=1.0\textwidth]{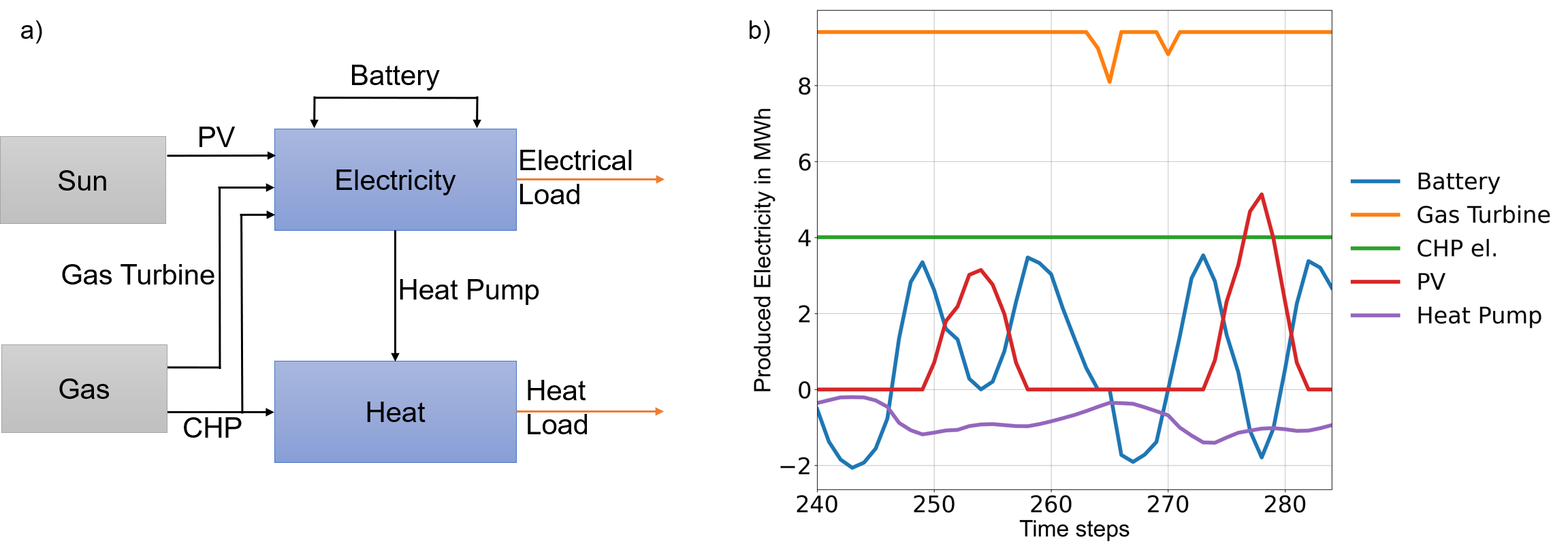}
    \caption{a) Exemplary Energy System b) Excerpt of the Plant Scheduling (electrical output or consumption)}
    \label{fig:002}
\end{figure}

\begin{table}[hbtp]
    \centering
    \begin{tabular}{lrrrrr}
    \hline
                                & \textbf{Battery} & \textbf{Gas Turbine} & \textbf{Combined heat and Power} & \textbf{Photovoltaic} & \textbf{Heat Pump} \\\hline
    P\textsubscript{inst.}in MW & 75               & 9.4                  & 4                                & 64.2                  & 4.3 \\\hline
    \end{tabular}
    \caption{Installed Capacities in the Optimised Energy System}
    \label{tab:001}
\end{table}

Linear programming is well suited for optimising those energy systems. Those problems feature a huge number of decision variables and constraints, but all constraints and the objective function can be expressed linearly regarding those decision variables. In energy system optimisation, the decision variables represent each component's scheduling, installed capacities, or other desired results. The constraints describe how single components could theoretically be deployed (e.g., their maximal load, their load change rate) and other boundary conditions like maximal CO2 emissions. The objective function typically is a cost function, which is minimised without violating the boundary conditions. After optimisation, the results show how to provide the desired loads in a cost-optimal way, i.e., the values of the decision variables representing the optimal scheduling and installed capacities. \autoref{fig:002} a) shows an exemplary system that aims to provide the electrical and heat load while minimising the total costs for PV, a gas turbine, a combined heat and power plant (CHP), a heat pump, and battery storage. \autoref{fig:002} b) shows the plant scheduling of two winter days with relatively low PV output. Table 1 shows the installed Capacities in this optimised energy system. This system is later used to present the equations as a practical example (see section Simple Energy System Optimisation Example).

% ----------------------------------------------------------------------------------------------------------------------------
% ----------------------------------------------------------------------------------------------------------------------------
% ----------------------------------------------------------------------------------------------------------------------------
% ----------------------------------------------------------------------------------------------------------------------------
% ----------------------------------------------------------------------------------------------------------------------------
\subsection{Mathematical Formulation of a Simple Problem}
In this study, the index c represents an element of the set of all components (not including storage, refer to Implementation of Storage) and the index t represents an element in the set of all time steps. Continuous decision variables are underlined (e.g. \underline{Pout\textsubscript{t,c}}) and decision variables that can only take on integer values are boxed (\fbox{on\textsubscript{t,c}}). Two decision variables are crucial. \underline{Pout\textsubscript{t,c}} is set up for each time step t and each component c and represents the respective load state. For example, this could represent the electrical output of a power plant or the hydrogen production of an electrolyser at time step t. If the capacity of component c also is optimised, another essential decision variable is \underline{Pinstalled\textsubscript{c}}, representing the installed capacity of the element. Depending on the problem, other decision variables are necessary. This could be the variable \fbox{on\textsubscript{t,c}} (representing whether component c is turned on at time step t) or \fbox{startup\textsubscript{t,c}} (representing whether component c was started between time step t and t-1). The chapter Equations using Mixed Integer Linear Programming explains the implementation of these variables. All decision variables are defined as non-negative.

Parameters are input values that have to be determined before the simulation. Those parameters become the coefficients in the linear optimisation problem. Examples are time series for the availability of wind and PV, efficiencies, or CO2 emission factors. A full description of the decision variables and parameters for the exemplary energy system can be found in the appendix.

After defining the decision variables, the constraints have to be formulated. \autoref{eq:001} forces that the load state of each component is always equal to or less than the available installed capacity. The installed capacity can be optimised, predetermined before the optimisation as an initial capacity, or a mix of those two options (as in \autoref{eq:001}). The installed capacity is multiplied by an availability factor. This equation can be set up for controllable producers like conventional power plants, fluctuating producers like PV, or wind turbines or any other component. For conventional power plants available\textsubscript{t,c} typically is a time series that is one if the plant is available at time step t and zero during maintenance. For PV and wind turbines, available\textsubscript{t,c} represents the specific availability of this energy source at a time step i.e. the ratio of produced electricity divided by the installed peak power. The less-than sign allows the optimiser to curtail the available production.
\begin{equation}\label{eq:001}
    \mathrm{Pout_{t,c} \leq \left( \underline{Pinstalled_{c}} + \underline{Pinstalled_{init,c}} \right) \cdot available_{t,c} \:\forall\: t \in time\:steps, \:\forall\: c \in components}
\end{equation}

\autoref{eq:002} ensures that a conservation balance at each node is fulfilled. The load at each time step has to equal everything that flows into the node minus everything that flows out of the node. The flow out of a node is calculated by dividing the respective components Pout\textsubscript{i,c} by its efficiency.
\begin{equation}\label{eq:002}
    \begin{gathered}
        \mathrm{Load_{t,n}=\sum_{c\,if\,component\,c\,starts\,at\,node\,n} \frac{1}{\eta_c} \cdot \underline{Pout_{t,c}} - \sum_{c\,if\,component\,c\,ends\,at\,node\,n} \underline{Pout_{t,c}}} \\
        \mathrm{\forall\: t \in time\:steps, \:\forall\: n \in nodes}
    \end{gathered}
\end{equation}

Finally, the cost function is expressed in \autoref{eq:003}. In this simple example, the investment costs c\textsubscript{inv}, maintenance costs c\textsubscript{maint} and fuel costs c\textsubscript{fuel} are considered. If the time span of the simulation is less than the lifetime of the components, costs like investment costs c\textsubscript{inv} have to be scaled down accordingly. Calculating the capital recovery factor with the per period interest rate i often is a suitable method (see \autoref{eq:004}).

Special attention must also be paid to the correct reference point of c\textsubscript{Inv}. In this study, the decision variables refer to the output of each component. Therefore, the costs need to be converted accordingly. If the costs are available for the input side, they have to be divided by the efficiency (see \autoref{eq:005}). This is typically the case for electrolysis, where costs are often given in $\frac{Euro}{kW_{el}}$ but are needed in $\frac{Euro}{kW_{H2}}$ for the simulation.
\begin{equation}\label{eq:003}
    \begin{gathered}
        \mathrm{c_{total} = \sum_{t\,in\,time\:steps;\,c\,in\,components} \underline{Pout_{t,c}} \cdot \frac{1}{\eta_c} \cdot \Delta t_t \cdot c_{fuel,t}} \\
        \mathrm{+ \sum_{c\,in\,components} \underline{Pinstalled_c} \cdot \left( c_{inv} + c_{maint} \right)}
    \end{gathered}
\end{equation}
\begin{equation}\label{eq:004}
    \mathrm{c_{Inv} = c_{Inv,total} \cdot \frac{\left( 1 + i \right)^n \cdot i}{\left( 1 + i \right)^n - 1}}
\end{equation}
\begin{equation}\label{eq:005}
    \mathrm{c_{Inv} = \frac{cinput_{Inv}}{\eta} + coutput_{Inv}}
\end{equation}

Those Equations describe a simple but complete energy system model that can be solved using a linear programming algorithm. 

% ----------------------------------------------------------------------------------------------------------------------------
% ----------------------------------------------------------------------------------------------------------------------------
% ----------------------------------------------------------------------------------------------------------------------------
% ----------------------------------------------------------------------------------------------------------------------------
% ----------------------------------------------------------------------------------------------------------------------------
\subsection{Advanced Equation for Linear Programming}
Although the model above already represents a complete system, further constraints can be implemented to represent a more realistic energy system.

% ----------------------------------------------------------------------------------------------------------------------------
\subsubsection{Capping the maximal installed capacity}
To cap the maximal installed capacity of any technology, \autoref{eq:006} can be implemented. This is usually applied when the potential of renewable technologies like PV or wind is limited.
\begin{equation}\label{eq:006}
    \mathrm{\underline{Pinstalled_c} \leq Pmaxinstalled_c \:\forall\:c \in components}
\end{equation}

% ----------------------------------------------------------------------------------------------------------------------------
\subsubsection{Multiple Input or Outputs of a Component}
\begin{figure}[hbtp]
    \centering
    \includegraphics[width=.90\textwidth]{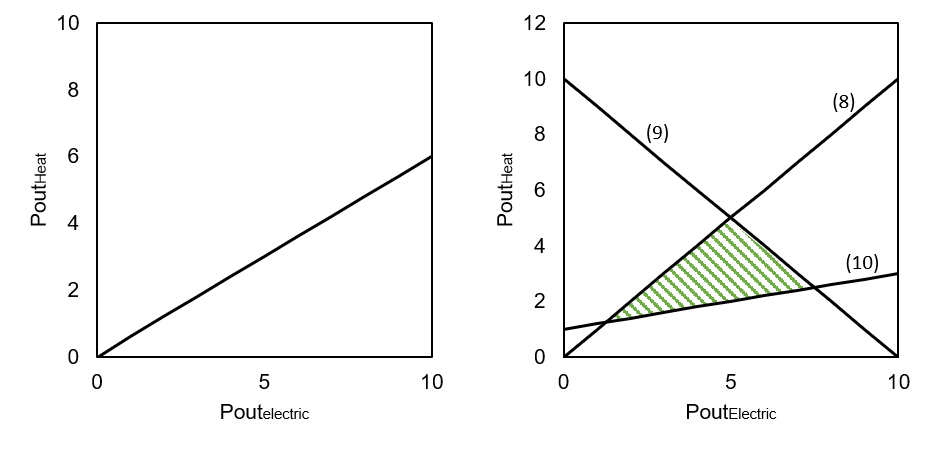}
    \caption{Characteristic Curve (a) and Characteristic Field (b)}
    \label{fig:003}
\end{figure}

Components may have multiple outputs, like CHP plants. There are two possibilities; either a characteristic curve or a characteristic field determines the dependency of those variables (see \autoref{fig:003}). To make the equations more readable, these types of equations are demonstrated using a combined heat and power plant as an example. Nevertheless, the idea of the equations can also be used for other components, more than two outputs, or multiple inputs.

If the link between electricity and heat is a line (see \autoref{fig:003} a)), the heat production can be expressed using the decision variable for the electricity production \underline{Pout\textsubscript{t,electric}} and the relevant efficiencies. Therefore, there is no need to define a dedicated variable for the heat output of this component. This term has to be integrated into the corresponding node balance (last term in \autoref{eq:007}).

\begin{equation}\label{eq:007}
    \begin{gathered}
        \mathrm{Load_{t,n} = \left[ \sum_{c\,if\,c\,starts\,at\,n} \frac{1}{\eta_c} \cdot \underline{Pout_{t,c}} - 
        \sum_{c\,if\,c\,ends\,at\,n} \underline{Pout_{t,c}}\right] + \sum_{c\,if\,c\,supplies\,heat\,to\,n} \frac{\eta_{heat}}{\eta_{electric}} \cdot \underline{Pout_{t,electric}}} \\
        \mathrm{\forall\:t \in time\:steps,\:\forall\:n \in nodes,\:\forall\:c \in components}
    \end{gathered}
\end{equation}

A new decision variable \underline{Pout\textsubscript{t,heat}} must be defined if there is a characteristic field (\autoref{fig:003} b)),. At least three inequalities define the field using both greater than and less than signs (e.g., \autoref{eq:008} to \ref{eq:010}), but there could be infinitely many. Nevertheless, it will always be convex since linear equations define the field.  The new decision variable \underline{Pout\textsubscript{t,heat}} has to be accounted for in the node balance (last term in \autoref{eq:011}).
\begin{equation}\label{eq:008}
    \mathrm{\underline{Pout_{t,heat}} \leq a_1 \cdot \underline{Pout_{t,electric}} + b_1\:\forall\:t\:\in\:timesteps}
\end{equation}
\begin{equation}\label{eq:009}
    \mathrm{\underline{Pout_{t,heat}} \leq a_2 \cdot \underline{Pout_{t,electric}} + b_2\:\forall\:t\:\in\:timesteps}
\end{equation}
\begin{equation}\label{eq:010}
    \mathrm{\underline{Pout_{t,heat}} \geq a_3 \cdot \underline{Pout_{t,electric}} + b_3\:\forall\:t\:\in\:timesteps}
\end{equation}
\begin{equation}\label{eq:011}
    \begin{gathered}
        \mathrm{Load_{t,n} = \left[ \sum_{c\,if\,c\,starts\,at\,n} \frac{1}{\eta_c} \cdot \underline{Pout_{t,c}} - 
        \sum_{c\,if\,c\,ends\,at\,n} \underline{Pout_{t,c}}\right] + \sum_{c\,if\,c\,supplies\,heat\,to\,n} \underline{Pout_{t,heat}}} \\
        \mathrm{\forall\:t \in time\:steps,\:\forall\:n \in nodes,\:\forall\:c \in components}
    \end{gathered}
\end{equation}

% ----------------------------------------------------------------------------------------------------------------------------
\subsubsection{Implementation of Storage}
Modelling energy storage can be done using \autoref{eq:012} to \ref{eq:015}. For storage, it is helpful to define two decision variables. One represents the flow into the storage \underline{Pcharge\textsubscript{i,s}} and one represents the flow out of the storage \underline{Pdischarge\textsubscript{i,s}}. A decision variable \underline{Capacity\textsubscript{s}} must be defined if the capacity also needs to be optimised. \autoref{eq:012} and \ref{eq:013} ensure, that the fill level is non-negative and less than its maximum capacity at any time step. These equations sum up the initial charge level as well as all charge- and discharge-flows that occurred until a time step. \autoref{eq:014} and \ref{eq:015} ensure that the maximum charge and discharge rates are not exceeded. The value for Capacity can either be set to a constant before the simulation or be defined as an optimisation variable. The values for Pmaxcharge\textsubscript{j} and Pmaxdischarge\textsubscript{j} can also be predefined, set to a fixed ratio of Capacity (\underline{Capacity}$=$C$\cdot$Pmaxcharge\textsubscript{j}; for example, if a battery has a given C-ratio), or can be set as independent optimisation variables.
\begin{equation}\label{eq:012}
    \begin{gathered}
        \mathrm{Wstored_{init,s} + \sum_{u\,in\,t} \left( \underline{Pcharge_{u,s}} \cdot \eta charge_s \cdot \Delta t_t - \underline{Pdischarge_{u,s}} \cdot \frac{1}{\eta discharge_s} \cdot \Delta t_t \right)} \\
        \mathrm{\geq 0} \\
        \mathrm{\forall\:t \in time\:steps,\:\forall\:s \in storages}
    \end{gathered}
\end{equation}
\begin{equation}\label{eq:013}
    \begin{gathered}
        \mathrm{Wstored_{init,s} + \sum_{u\,in\,t} \left( \underline{Pcharge_{u,s}} \cdot \eta charge_s \cdot \Delta t_t - \underline{Pdischarge_{u,s}} \cdot \frac{1}{\eta discharge_s} \cdot \Delta t_t \right)} \\
        \mathrm{\geq \underline{Capacity_s} + Capacity_s} \\
        \mathrm{\forall\:t \in time\:steps,\:\forall\:s \in storages}
    \end{gathered}
\end{equation}
\begin{equation}\label{eq:014}
    \begin{gathered}
        \mathrm{\underline{Pcharge_{t,s}} \leq Pmaxcharge_s + \frac{1}{C_s} \cdot \underline{Capacity_s} + \underline{Pmaxcharge_s}} \\
        \mathrm{\forall\:t \in time\:steps,\:\forall\:s \in storages}
    \end{gathered}
\end{equation}
\begin{equation}\label{eq:015}
    \begin{gathered}
        \mathrm{\underline{Pdischarge_{t,s}} \leq Pmaxdischarge_s + \frac{1}{C_s} \cdot \underline{Capacity_s} + \underline{Pmaxdischarge_s}} \\
        \mathrm{\forall\:t \in time\:steps,\:\forall\:s \in storages}
    \end{gathered}
\end{equation}

% ----------------------------------------------------------------------------------------------------------------------------
\subsubsection{Load Change Rate}
To limit the load change rate of components like power plants, \autoref{eq:016} to \ref{eq:017} can be implemented. In this case, it is essential to understand that in those simulations, the load change resembles the maximum difference of the averages of the output of two successive time steps. The maximum load change usually is expressed as a fixed percentage of the decision variable \underline{Pinstalled\textsubscript{j}}. For example, a power plant might have a load change rate of lcrup$=80\,\frac{\%}{h}$. Nevertheless, it would be conceivable to optimise how fast a power plant can ramp up or down if the goal is to know how flexible future power plants have to be (e.g., using a decision variable \underline{LCRup}).
\begin{equation}\label{eq:016}
    \begin{gathered}
        \mathrm{Pout_{t,c} - Pout_{t-1,c} \leq lcrup_c \cdot \underline{Pinstalled_j} + \underline{LCRup_c}} \\
        \mathrm{\forall\:t \in time\:steps:t>1,\:\forall\:c \in components}
    \end{gathered}
\end{equation}
\begin{equation}\label{eq:017}
    \begin{gathered}
        \mathrm{Pout_{t-1,c} - Pout_{t,c} \leq lcrdown_c \cdot \underline{Pinstalled_j} + \underline{LCRdown_c}} \\
        \mathrm{\forall\:t \in time\:steps:t>1,\:\forall\:c \in components}
    \end{gathered}
\end{equation}

% ----------------------------------------------------------------------------------------------------------------------------
\subsubsection{Multiple Building Periods}
For larger simulations, it is also possible to define multiple building periods p, for example, to investigate in which years specific technologies must be implemented into the energy system. Then, the decision variable \underline{Pinstalled\textsubscript{p,c}} has to be defined for each building period p and each component. \autoref{eq:018} has to be set up if the time step t falls into the building period p. This equation replaces \autoref{eq:001}.
To obtain a variable representing the amount of capacity added from one building period to the next, \autoref{eq:019} can be set up. Two mechanisms combine for this to work: First, \underline{Pbuilt\textsubscript{p,c}} needs to be associated with positive costs in the objective function so that the optimiser will always choose the lowest possible value. Second, this lowest possible value is defined by \autoref{eq:019} or the non-negative condition of \underline{Pbuilt\textsubscript{p,c}}. By combining these mechanisms, \underline{Pbuilt\textsubscript{p,c}} equals \underline{Pinstalled\textsubscript{p,c}}$-$\underline{Pinstalled\textsubscript{p-1,c}} if capacity is built. Otherwise, \underline{Pbuilt\textsubscript{p,c}} is zero. 
\begin{equation}\label{eq:018}
    \mathrm{\underline{Pout_{t,c}} \leq \underline{Pinstalled_{p,c}} \cdot available_{t,c}}
\end{equation}
\begin{equation}\label{eq:019}
    \begin{gathered}
        \mathrm{\underline{Pinstalled_{p,c}} - \underline{Pinstalled_{p-1, c}} \leq \underline{Pbuilt_{p, c}}}\\
        \mathrm{\forall\:p \in periods:p>1,\:\forall\:c \in components}
    \end{gathered}
\end{equation}

% ----------------------------------------------------------------------------------------------------------------------------
\subsubsection{Cost Function}
The cost function must also account for all those decision variables (see \autoref{eq:020}). If this were not the case, the optimizer could set decision variables like the maximum capacity of storages to a very high value regardless of the problem. The results would, therefore, not yield any additional insights into the system.
\begin{equation}\label{eq:020}
    \begin{gathered}
        \mathrm{C_{total} = ...}\\
        \mathrm{+ c_{capacity,j} \cdot \underline{Capacity_j}} \\
        \mathrm{+ c_{Pmaxcharge,j} \cdot \underline{Pmaxcharge_j}}\\
        \mathrm{+ c_{maxLoadChangeRatedown} \cdot \underline{LCRup_c}} \\
        \mathrm{+ ...}
    \end{gathered}
\end{equation}

% ----------------------------------------------------------------------------------------------------------------------------
\subsubsection{Additional Constraints}
Further constraints not concerning specific technologies can be imposed on the model. However, they are often very problem-specific, and there is a multitude of possibilities. Therefore, only one option is shown here: imposing a limit on total CO2 emissions (\autoref{eq:021}). A different way to punish CO2 emissions would be to include those emissions in the cost function.
\begin{equation}\label{eq:021}
    \mathrm{\sum_{t\,in\,time\,steps;\,c\,in\,components} \underline{Pout_{t,c}} \cdot \frac{1}{\eta_c} \cdot specCO_2Emissions_c \leq CO_2max}
\end{equation}

% ----------------------------------------------------------------------------------------------------------------------------
% ----------------------------------------------------------------------------------------------------------------------------
% ----------------------------------------------------------------------------------------------------------------------------
% ----------------------------------------------------------------------------------------------------------------------------
% ----------------------------------------------------------------------------------------------------------------------------
\subsection{Equations using Mixed Integer Linear Programming}
Mixed Integer Linear Programming allows the introduction of integer decision variables, variables that can only take on whole numbers. This choice opens up the possibility of implementing discrete steps.

% ----------------------------------------------------------------------------------------------------------------------------
\subsubsection{Discrete Startups}
Using \autoref{eq:022} to \ref{eq:024} allows modelling a power plant with discrete startups when the variable \fbox{on\textsubscript{t,c}} is defined as an integer decision variable (usually binary). \autoref{eq:022} and \ref{eq:023} force \underline{Pout\textsubscript{t,c}} to be zero if \fbox{on\textsubscript{t,c}} is zero (if the power plant is turned off). Those equations also force \underline{Pout\textsubscript{t,c}} to be larger than the minimum load of a single unit but less than the installed capacity of a single unit if exactly one power plant is turned on (\fbox{on\textsubscript{t,c}}$=$1). However, the equations also work if multiple units are represented. Then, \underline{Pout\textsubscript{t,c}} needs to be larger than the combined minimum load of the number of units that are turned on but less than their combined installed capacity. In \autoref{eq:022} and \ref{eq:023}, Pinstalledunit\textsubscript{c} needs to be a predetermined factor to keep the equation linear regarding the decision variables. However, optimising the number of plants (see \autoref{eq:027}) and therefore the total installed capacity is possible.

\autoref{eq:024} defines a decision variable \fbox{startup\textsubscript{t,c}}. This equation works similarly to \autoref{eq:029} The variable \fbox{startup\textsubscript{t,c}} needs to be associated with positive costs in the objective function. Then the optimiser will choose the minimum valid value. The minimum value corresponds to the number of units turned on if plants are started (if \fbox{on\textsubscript{t,c}}$-$\fbox{on\textsubscript{t-1,c}} is positive) and zero due to the non-negative condition if no plants are turned on (then \fbox{on\textsubscript{t,c}}$-$\fbox{on\textsubscript{t-1,c}} is zero or negative). If implemented, these equations replace \autoref{eq:001}.

If the node balance is adjusted according to \autoref{eq:025}, a partial load efficiency can be implemented with the factors p\textsubscript{c} and o\textsubscript{c}. This enables the implementation of plants with lower efficiency at partial load. The efficiency is then dependent on the load (see \autoref{eq:026})

\begin{equation}\label{eq:022}
    \begin{gathered}
        \mathrm{\underline{Pout_{t,c}} \leq \boxed{\mathrm{on_{t,c}}} \cdot Pinstalledunit_c \cdot available_{t,c}} \\
        \mathrm{\forall\:t \in time\:steps,\:\forall\:c \in components}
    \end{gathered}
\end{equation}
\begin{equation}\label{eq:023}
    \begin{gathered}
        \mathrm{\underline{Pout_{t,c}} \geq \boxed{\mathrm{on_{t,c}}} \cdot Pminunit_c} \\
        \mathrm{\forall\:t \in time\:steps,\:\forall\:c \in components}
    \end{gathered}
\end{equation}
\begin{equation}\label{eq:024}
    \begin{gathered}
        \mathrm{\boxed{\mathrm{on_{t,c}}} - \boxed{\mathrm{on_{t-1,c}}} \leq \boxed{\mathrm{startup_{t,c}}}} \\
        \mathrm{\forall\:t \in time\:steps,\:\forall\:c \in components}
    \end{gathered}
\end{equation}
\begin{equation}\label{eq:025}
    \begin{gathered}
        \mathrm{Load_{t,n} = \sum_{c\,if\,c\,starts\,at\,n} \left( p_c \cdot \underline{Pout_{t,c}} + o_c \cdot \boxed{\mathrm{on_{t,c}}} \right) - \sum_{c\,if\,c\,ends\,at\,n} \underline{Pout_{t,c}}} \\
        \mathrm{\forall\:t \in time\:steps,\:\forall\:n \in nodes}
    \end{gathered}
\end{equation}
\begin{equation}\label{eq:026}
    \begin{gathered}
        \eta_{t,c} = \frac{\underline{Pout_{t,c}}}{p_c \cdot \underline{Pout_t,c} + o_c}
    \end{gathered}
\end{equation}

% ----------------------------------------------------------------------------------------------------------------------------
\subsubsection{Discrete Number of units}
If units are supposed to be added discretely, \autoref{eq:027} is added. This indirectly allows optimising the installed capacity by optimising the number of units.
\begin{equation}\label{eq:027}
    \mathrm{\boxed{\mathrm{on_{t,c}}} \leq \boxed{\mathrm{units_{c}}} \: \forall\:t \in time\:steps,\:\forall\:c \in components}
\end{equation}

% ----------------------------------------------------------------------------------------------------------------------------
\subsubsection{Minimum Down-time and Minimum Up-time} 
If the variable \fbox{on\textsubscript{t,c}} is not only integer but binary, constraints for minimum downtime and minimum up-time of N time steps can be introduced. These equations will be explained using minimum downtime \autoref{eq:028}, however, minimum up-time works analogously \autoref{eq:029}. If the plant is not turned on before time step i, the equation does not have an effect. Then, the left-hand side of the equation is zero, since \fbox{on\textsubscript{t,c}} $=$ \fbox{on\textsubscript{t-1,c}} and the right-hand side will always be larger than or equal to zero. If the plant is turned on before time step I, \fbox{on\textsubscript{t}} $=1$ and \fbox{on\textsubscript{t-1}} $=0$ and the left-hand side will therefore be N. Thus, this constraint will only not be violated if $\mathrm{\sum_{m\,in\,n} on_{t-m}}$ is zero. Then, the plant was turned off for the last N time steps and the minimum downtime was fulfilled. 
\begin{equation}\label{eq:028}
    \begin{gathered}
        \mathrm{\left( \boxed{\mathrm{on_{t,c}}} - \boxed{\mathrm{on_{t-1,c}}} \right) \cdot N \leq N -\sum_{m\,in\,N} \boxed{\mathrm{on_{t-m,c}}} } \\
        \mathrm{\forall\:t \in time\:steps,\:\forall\:c \in components}
    \end{gathered}
\end{equation}
\begin{equation}\label{eq:029}
    \begin{gathered}
        \mathrm{\left( \boxed{\mathrm{on_{t-1,c}}} - \boxed{\mathrm{on_{t,c}}} \right) \cdot N \leq \sum_{m\,in\,N} \boxed{\mathrm{on_{t-m,c}}} } \\
        \mathrm{\forall\:t \in time\:steps,\:\forall\:c \in components}
    \end{gathered}
\end{equation}

% ----------------------------------------------------------------------------------------------------------------------------
% ----------------------------------------------------------------------------------------------------------------------------
% ----------------------------------------------------------------------------------------------------------------------------
% ----------------------------------------------------------------------------------------------------------------------------
% ----------------------------------------------------------------------------------------------------------------------------
\subsection{Simple Energy System Optimisation Example}
The system shown in \autoref{fig:002} a) represents a small exemplary energy system. The purpose of this example is to show the mathematical formulation of the optimization problem. Real energy systems are usually much more complex. The model aims to determine how the required electricity and heat can be provided cost-efficiently using the available technologies. The regarded year is divided into 8760 time steps, each representing one hour. Relevant information regarding the nodes and edges can be found in the appendix. Hourly load balances are calculated in MW for both nodes in each time step. The system is allowed to install five different components:

\begin{itemize}
    \item A PV plant producing electricity with a given hourly load factor
    \item A combined heat and power plant, burning natural gas to produce electricity as well as heat
    \item A gas turbine burning natural gas to produce electricity
    \item A battery to store electricity 
    \item A heat pump using electricity to produce heat. 
\end{itemize}

As a result, the installed capacities and the scheduling are calculated. \autoref{fig:002} b) depicts an extract of electricity production of 48 hours. \autoref{tab:001} Table 1 shows the resulting installed capacities. In a cost-optimal system, the gas turbine and the CHP are operated in base-load to provide heat or/and electricity, respectively. Fluctuations in the heat demand are matched via the heat pump. Electricity provided by the gas turbine and the CHP is complemented by a combination of PV and battery storage. The operation of the battery supplements the electricity production from PV.

As a summary, the full set of equations is shown in \autoref{eq:030} (objective function) and \autoref{eq:031}:
\begin{equation}\label{eq:030}
    \begin{gathered}
        \mathrm{min \Biggl[ \sum_{t\,in\,1..8760} \left( \frac{1}{\eta CHP} \cdot \underline{P_{8760,CHP}} \cdot \Delta t_t \cdot c_{gas} + \frac{1}{\eta Gas\,Turbine} \cdot \underline{P_{8760,Gas\,Turbine}} \cdot \Delta t_t \cdot c_{gas} \right)} \\
        \mathrm{+ \underline{Pinstalled_{PV}} \cdot c_{Inv,PV} + \underline{Pinstalled_{CHP}} \cdot c_{Inv,CHP}} \\
        \mathrm{+ \underline{Pinstalled_{Gas\,Turbine}} \cdot c_{Inv,Gas\,Turbine} + \underline{Pinstalled_{Heat\,Pump}} \cdot c_{Inv,Heat\,Pump}} \\
        \mathrm{+ \underline{Pinstalled_{Battery}} \cdot c_{Inv,Battery} \Biggr]}
    \end{gathered}
\end{equation}

% ----------------------------------------------------------------------------------------------------------------------------
% ----------------------------------------------------------------------------------------------------------------------------
% ----------------------------------------------------------------------------------------------------------------------------
% ----------------------------------------------------------------------------------------------------------------------------
% ----------------------------------------------------------------------------------------------------------------------------
\section{Conclusion}
This paper provides a detailed explanation of the methodology of energy system simulation. Moreover, the implantation of a simplified energy system is shown. Therefore, this study is a valuable tool for researchers trying to understand the methodology of energy system simulations and building their framework.

First, this paper explains the representation of the energy system as nodes and edges and the basic transfer of energy or mass between different nodes. Then, the primary constraints are defined, advanced tools are introduced, and possibilities using mixed integer linear programming are shown. Lastly, the equations are shown for an exemplary system to show them from a more practical point of view.

%\printbibliography

\bibliographystyle{unsrtnat}
%\bibliography{references}  %%% Uncomment this line and comment out the ``thebibliography'' section below to use the external .bib file (using bibtex) .


\begin{thebibliography}{1}

\bibitem{Erlach.2018}
Erlach, Berit and Henning, Hans-Martin and Kost, Christoph and Palzer, Andreas and Stephanos, Cyril.
acatech.
\textit{Optimierungsmodell REMod-D: Materialien zur Analyse <<Sektorkopplung>> -- Untersuchungen und {\"U}berlegungen zur Entwicklung eines integrierten Energiesystems}.
Schriftenreihe Energiesysteme der Zukunft. 2018.

\bibitem{Hughes.2021}
Hughes, Paul.
{International Energy Agency}.
\textit{World Energy Model Documentation: 2020 Version}. 2021.

\bibitem{Loulou.2016}
Loulou, Richard and Goldstein, Gary and Kanudia, Amit and Lettila, Antti and Remme, Uwe.
{International Energy Agency}.
\textit{Documentation for the TIMES Model}. 2016.

\bibitem{Miehling.2021}
Miehling, Sebastian and Schweiger, Benedikt and Wedel, Wolf and Hanel, Andreas and Schweiger, Jakob and Schwermer, Rene and Blume, Maximilian and Spliethoff, Hartmut.
\textit{100 {\%} erneuerbare Energien f{\"u}r Bayern: Potenziale und Strukturen einer Vollversorgung in den Sektoren Strom, W{\"a}rme und Mobilit{\"a}t}. 2021

\bibitem{Miehling.2022}
Miehling, Sebastian and Fendt, Sebastian and Spliethoff, Hartmut.
\textit{Optimal integration of Power-to-X plants in a future European energy system and the resulting dynamic requirements}.
In Energy Conversion and Management, \textit{doi:} 10.1016/j.enconman.2021.115020

\end{thebibliography}

\clearpage
\appendix
\section{Appendix}
\begin{table}[hbtp]
    \caption{Parameter overview of the considered technologies in the energy system example.}
    \label{tab:002}
    \begin{tabular}{lrrrrr}\hline
                                           & \multicolumn{4}{l}{\textbf{Linear Element}s}                               & \textbf{Storage}     \\
                                           &\textbf{Heat Pump} & \textbf{Gas Turbine}     & \textbf{CHP}         & \textbf{PV}              & \textbf{Battery}     \\\hline
    Starting Node                          & Electricity     & Gas             & Gas         & Sun             & Electricity \\
    Ending Node 1                          & Heat            & Electricity     & Electricity & Electricity     & Electricity \\
    Ending Node 2                          & /               & /               & Heat        & /               & /           \\
    Maximal Capacity                       & 1000 MW         & 1000 MW         & 1000 MW     & 1000 MW         & 1000 MW     \\
    Minimal Capacity                       & 0 MW            & 0 MW            & 0 MW        & 0 MW            & 0 MW        \\
    Installed Capacity                     & 0 MW            & 0 MW            & 0 MW        & 0 MW            & 0 MW        \\
    Max. Output/Inst. Capacity             & 1 MW            & 1 MW            & 1 MW        & 1 MW            & /           \\
    Min. Output/Inst. Capacity             & 0 MW            & 0 MW            & 0 MW        & 0 MW            & /           \\
    Charge Flow                            & /               & /               & /           & /               & 1 MW        \\
    Discharge Flow                         & /               & /               & /           & /               & 1 MW        \\
    Load Change Rate                       & 1               & 1               & 1           & /               & /           \\
    Efficiency 1                           & 3               & 0,4             & 0,37        & 1               & 0,9604      \\
    Efficiency 2                           & /               & /               & 0,48        & /               & /           \\
    Ratio Node 1/Node2                     & /               & /               & 0,768       & /               & /           \\
    Fuel Costs in €/MWh                    & 21,61           & 21,61           & 21,61       & 0               & /           \\
    CO2   equivalent in kg/MWh             & 0,202           & 0,202           & 0,202       & 0               & /           \\
    CO2-Price in €/kg                      & 30              & 30              & 30          & /               & /           \\
    Availability                           & 1               & 1               & 1           & Availability PV & 1           \\
    Annual Inv. Costs in €/MWa             & 19028           & 24850           & 45795       & 21300           & 8520        \\\hline
    \end{tabular}
\end{table}

\begin{equation}\label{eq:031}
    \begin{gathered}
        \mathrm{\underline{P_{1,PV}} \leq \underline{Pinstalled_{PV}} \cdot available_{1,PV}} \\
        \mathrm{\underline{P_{1,CHP}} \leq \underline{Pinstalled_{CHP}}} \\
        \mathrm{\underline{P_{1,Gas\,Turbine}} \leq \underline{Pinstalled_{Gas\,Turbine}}} \\
        \mathrm{\underline{P_{1,Heat\,Pump}} \leq \underline{Pinstalled_{Heat\,Pump}}} \\
        \mathrm{\underline{Pcharge_{1,Battery}} \leq \frac{1}{C_{Battery}} \underline{Capacity_{Battery}}} \\
        \mathrm{\underline{Pdischarge_{1,Battery}} \leq \frac{1}{C_{Battery}} \underline{Capacity_{Battery}}} \\
        \mathrm{\sum_{u\,in\,1..1} \left( \underline{Pcharge_{u,s}} \cdot \eta charge_s \cdot \Delta t_t - \underline{Pdischarge_{u,s}} \cdot \frac{1}{\eta discharge_s} \cdot \Delta t_t \right) \geq 0} \\
        \mathrm{\sum_{u\,in\,1..1} \left( \underline{Pcharge_{u,s}} \cdot \eta charge_s \cdot \Delta t_t - \underline{Pdischarge_{u,s}} \cdot \frac{1}{\eta discharge_s} \cdot \Delta t_t \right) \geq \underline{Capacity_s}} \\
        \mathrm{Load_{Elec,1} = \underline{P_{1,PV}} + \underline{P_{1,CHP}} + \underline{P_{1,Gas\,Turbine}} - \frac{1}{\eta_{Heat\,Pump}} \cdot \underline{P_{1,Heat\,Pump}} + Pdischarge_{1,Battery}} \\
        \mathrm{- Pcharge_{1,Battery}} \\
        \mathrm{Load_{Heat,1} = \frac{1}{\eta_{Heat\,Pump}} \cdot \underline{P_{1,Heat\,Pump}} + \frac{\eta_{Heat,CHP}}{\eta_{Elec,CHP}} \cdot \underline{P_{1,CHP}}}\\
        ...\\
        \mathrm{\underline{P_{8760,PV}} \leq \underline{Pinstalled_{PV}} \cdot available_{8760,PV}} \\
        \mathrm{\underline{P_{8760,CHP}} \leq \underline{Pinstalled_{CHP}}} \\
        \mathrm{\underline{P_{8760,Gas\,Turbine}} \leq \underline{Pinstalled_{Gas\,Turbine}}} \\
        \mathrm{\underline{P_{8760,Heat\,Pump}} \leq \underline{Pinstalled_{Heat\,Pump}}} \\
        \mathrm{\underline{Pcharge_{8760,Battery}} \leq \frac{1}{C_{Battery}} \underline{Capacity_{Battery}}} \\
        \mathrm{\underline{Pdischarge_{8760,Battery}} \leq \frac{1}{C_{Battery}} \underline{Capacity_{Battery}}} \\
        \mathrm{\sum_{u\,in\,1..8760} \left( \underline{Pcharge_{u,s}} \cdot \eta charge_s \cdot \Delta t_t - \underline{Pdischarge_{u,s}} \cdot \frac{1}{\eta discharge_s} \cdot \Delta t_t \right) \geq 0} \\
        \mathrm{\sum_{u\,in\,1..8760} \left( \underline{Pcharge_{u,s}} \cdot \eta charge_s \cdot \Delta t_t - \underline{Pdischarge_{u,s}} \cdot \frac{1}{\eta discharge_s} \cdot \Delta t_t \right) \geq \underline{Capacity_s}} \\
        \mathrm{Load_{Elec,8760} = \underline{P_{8760,PV}} + \underline{P_{8760,CHP}} + \underline{P_{8760,Gas\,Turbine}} - \frac{1}{\eta_{Heat\,Pump}} \cdot \underline{P_{8760,Heat\,Pump}}} \\
        \mathrm{+ Pdischarge_{8760,Battery} - Pcharge_{8760,Battery}} \\
        \mathrm{Load_{Heat,8760} = \frac{1}{\eta_{Heat\,Pump}} \cdot \underline{P_{8760,Heat\,Pump}} + \frac{\eta_{Heat,CHP}}{\eta_{Elec,CHP}} \cdot \underline{P_{8760,CHP}}}\\
        ...\\
        \mathrm{(see\:section\:Advanced\:Equation\:for\:Linear\:Programming)}
    \end{gathered}
\end{equation}

\begin{landscape}
\begin{table}[hbtp]
    \caption{Nomenclature (1/2)}
    \label{tab:003}
    \begin{tabular}{lllp{13cm}}\hline
    & \textbf{Symbol} & \textbf{Name} & \textbf{Description} \\ \hline
    % -------------------------------------------------------------------------------------------------------------------------------------------------
    \multirow{12}{*}{\rotatebox[origin=c]{90}{\textbf{Decision Variable}}}
    & \underline{Capacity}      & Storage Capacity & Capacity of a storage \\ 
    & \underline{LCRdown}       & Load Change Rate Down & Variable that indicates, how much the power output of a component can decrease from one time step to the next one \\ 
    & \underline{LCRup}         & Load Change Rate Up & Variable that indicates, how much the power output of a component can increase from one time step to the next one \\ 
    & \fbox{on}                 & On & Variable that represents, whether a component is turned on \\ 
    & \underline{Pcharge}       & Power Charging & Number of power that a storage is charged with \\ 
    & \underline{Pdischarge}    & Power Discharging & Number of power that a storage is discharged with \\ 
    & \underline{Pinstalled}    & Installed Capacity of a Component & Installed Capacity of a component. \\ 
    & \underline{Pmaxcharge}    & Max. Charging Power & Max. charging power of a storage \\ 
    & \underline{Pmaxdischarge} & Max. Discharging Power & Max. discharging power of a storage \\ 
    & \underline{Pout}          & Power Output & The stream that a component provides to a node. \\ 
    & \fbox{startup}            & Startup & Variable that represents, whether a component was started between time step t and t-1 \\ 
    & \fbox{units}              & Number of Installed Units & Number of installed units \\ \hline
    % -------------------------------------------------------------------------------------------------------------------------------------------------
    \multirow{6}{*}{\rotatebox[origin=c]{90}{\textbf{Index}}}
    & c                         & Item in the Set of all Components & \\
    & init                      & Initial & Marks initial values: the values before the first time step \\
    & n                         & Item in the Set of all Nodes & \\
    & p                         & Item in the Set of all Building periods & \\
    & s                         & Item in the Set of all Storages & \\
    & t                         & Item in the Set of all time Steps & \\ \hline
    % -------------------------------------------------------------------------------------------------------------------------------------------------
    \multirow{7}{*}{\rotatebox[origin=c]{90}{\textbf{Parameter}}}
    & $\Delta$t & Time Step Size & Step size of a time step (e.g. 1 hour or 15 minutes) \\
    & $\eta$charge & Charging Efficiency & Efficiency of charging a storage \\
    & $\eta$discharge & Discharging Efficiency & Efficiency of discharging a storage \\
    & a &  & Parameter to design the characteristic field \\
    & available & Available & Represents the share of the installed capacity that is available at a specific time step \\
    & b &  & Parameter to design the characteristic field \\
    & C\textsubscript{total} & Total Costs & Total Costs \\ \hline
    \end{tabular}
\end{table}    
\end{landscape}

\begin{landscape}
\begin{table}[hbtp]
    \caption{Nomenclature (2/2)}
    \label{tab:004}
    \begin{tabular}{lllp{13cm}}\hline
    & \textbf{Symbol} & \textbf{Name} & \textbf{Description} \\ \hline    
    \multirow{18}{*}{\rotatebox[origin=c]{90}{\textbf{Parameter}}}

    & c\textsubscript{x} & Specific Costs & Specific costs of a decision variable \\
    & Capacity & Storage Capacity & Capacity of a storage \\
    & CO\textsubscript{2} max & Max. CO2 Emissions & Max. CO2 Emissions that are allowed in a simulation \\
    & $\eta$ & Efficiency & Efficiency of a component i.e. the ratio of the provided stream and the consumed stream \\
    & lcrdown & Specific Load Change Rate Down & Determines how much the power output of a component can decrease from one time step to the next one (relative to the installed capacity) \\
    & lcrup & Specific Load Change Rate Up & Determines how much the power output of a component can increase from one time step to the next one (relative to the installed capacity) \\
    & Load & Load & Time series of the load at a node \\
    & N &  & Number of time steps for min downtime and min up-time \\
    & o &  & Parameter to calculate the input power dependent on the output power in mixed integer linear programming \\
    & p &  & Parameter to calculate the input power dependent on the output power in mixed integer linear programming \\
    & Pinstalled & Installed Capacity of a Component & Installed Capacity of a component. This variable is a value predetermined by the user and represents the initial installed capacity. \\
    & Pinstalledunit & Installed Capacity of a Single Unit & capacity of a single unit in mixed integer programming \\
    & Pmaxcharge & Max. Charging Rate & Max. charging power of a storage \\
    & Pmaxdischarge & Max. Discharging Rate & Max. discharging power of a storage \\
    & Pmaxinstalled & Max. Installed Capacity & Max. installed capacity \\
    & Pminunit & Min. Load of a Single Unit & Min. load of a single unit in mixed integer programming \\
    & specCO\textsubscript{2} Emissions & Specific CO2 Emissions & Specific CO2 Emissions of a fuel \\ \hline
    \end{tabular}
\end{table}
\end{landscape}

\end{document}